\documentclass[11pt]{article}

\usepackage[margin=1.07in]{geometry}
\usepackage{amsmath,amssymb,amsthm,mathtools}
\usepackage{microtype}
\usepackage{enumitem}
\usepackage[hidelinks]{hyperref}

\newtheorem{theorem}{Theorem}[section]
\newtheorem{proposition}[theorem]{Proposition}
\newtheorem{lemma}[theorem]{Lemma}
\newtheorem{corollary}[theorem]{Corollary}
\theoremstyle{definition}
\newtheorem{definition}[theorem]{Definition}
\theoremstyle{remark}

\newcommand{\C}{\mathbb C}
\newcommand{\R}{\mathbb R}

\newcommand{\caplog}{\operatorname{cap}}
\newcommand{\Ree}{\operatorname{Re}}

\title{Sharpness of the Brjuno Condition for Hyperbolic Bishop Surfaces%
}

\author{Dijia Chen%
\thanks{Department of Mathematics,
University of Wisconsin--Madison,
Madison, WI 53706, USA.
E-mail: \texttt{dchen426@wisc.edu}.}}

\date{}

\begin{document}
\maketitle
\begin{abstract}
Let $\mu \in \mathbb{C}$ satisfy $|\mu|=1$, and assume that $\mu$ is not a root of unity and is not Brjuno. Let $Q_\gamma$ denote the Bishop quadric with $\gamma = 1/\sqrt{2+\mu+\mu^{-1}} > 1/2$. Gong and Stolovitch conjectured that there exists a real-analytic surface in $\mathbb{C}^2$ with Moser--Webster invariant $\mu$ that is formally but not holomorphically equivalent to $Q_\gamma$. We prove this conjecture.
\end{abstract}

\section{Introduction}

Let $Q_\gamma$ be the Bishop quadric in $\mathbb{C}^2$ defined in \cite{Bis65} by $z_2=z_1\overline{z}_1+\gamma(z_1^2+\overline{z}_1^{\,2})$, with $\gamma>1/2$. At the origin, the real tangent plane of $Q_\gamma$ is a complex line, so the origin is a complex tangent point. In general, Bishop classified such complex tangencies of $Q_\gamma$ and their higher order perturbations according to the invariant $\gamma$: elliptic when $0\leq\gamma<1/2$, parabolic when $\gamma=1/2$, and hyperbolic when $\gamma>1/2$. In the elliptic case, Bishop showed that higher-order perturbations of an elliptic Bishop quadric still admit attached holomorphic discs. These discs form families shrinking to the complex tangent point \cite{Bis65}. By contrast, a surface near a hyperbolic complex tangency is locally polynomially convex \cite{FS91}. In the parabolic case, however, local polynomial convexity may hold or fail, depending on the higher-order terms \cite{Jor97}. The reader should consult Kenig and Webster \cite{KW82} and Huang and Krantz \cite{HK95} for further developments of Bishop's work on the regularity of the family of attached analytic discs and local hulls of holomorphy of elliptic smooth or real analytic real surfaces in $\C^2$. For higher dimensional elliptic manifolds, see Kenig-Webster~\cite{KW84} for the finite smooth case and Huang \cite{Hua98} for a complete theory on the (finite) smooth and real analytic cases.

Beyond these geometric properties, another central problem concerns the formal and holomorphic equivalence of real-analytic surfaces near complex tangencies. Moser and Webster showed that these equivalence problems can be reduced to conjugacy problems for triples of involutions \cite{MW83}. For a real-analytic surface $M$ with a complex tangent of nonzero Bishop invariant, the associated triple is $(\tau_{1,M},\tau_{2,M},\rho_M)$. Here $\tau_{1,M}$ and $\tau_{2,M}$ are holomorphic involutions, while $\rho_M$ is an antiholomorphic involution encoding the real structure. The two holomorphic involutions arise from the natural projections of the complexification of $M$. These involutions satisfy $\tau_{2,M}=\rho_M\circ\tau_{1,M}\circ\rho_M$. Under the Moser--Webster correspondence, formal equivalence corresponds to formal simultaneous conjugacy of the associated triples. Holomorphic equivalence corresponds to holomorphic simultaneous conjugacy.

Let $\sigma_M=\tau_{1,M}\circ\tau_{2,M}$. Define $\sigma_{Q_\gamma}$ similarly for the Bishop quadric $Q_\gamma$. If $M$ is formally equivalent to $Q_\gamma$, then $\sigma_M$ and $\sigma_{Q_\gamma}$ are formally conjugate. If $M$ is holomorphically equivalent to $Q_\gamma$, then these two maps are holomorphically conjugate. Since $\sigma_{Q_\gamma}$ is linear in suitable coordinates, holomorphic equivalence to $Q_\gamma$ requires $\sigma_M$ to be holomorphically linearizable.

For a hyperbolic Bishop invariant $\gamma$, choose $\lambda\in\mathbb{C}$ with $|\lambda|=1$ and $\lambda+\overline{\lambda}=\gamma^{-1}$. Set $\mu=\lambda^2$. The linear part of $\sigma_M$ has eigenvalues $\mu$ and $\mu^{-1}$, with $\mu+\mu^{-1}=\gamma^{-2}-2$. Following Gong and Stolovitch \cite{GS26}, we call $\mu$ the Moser--Webster invariant. In coordinates where $\sigma_{Q_\gamma}$ is linear, the formal conjugacy equation is $H\circ\sigma_M=\sigma_{Q_\gamma}\circ H$. Even if this equation admits a formal solution, that solution need not converge. Indeed, when $\mu$ is not a root of unity, the coefficient equations may involve arbitrarily small divisors. This raises a natural question: does the formal equivalence of a real analytic surface with $Q_\gamma$ imply their holomorphic equivalence?

In 2004, Gong constructed a real-analytic hyperbolic Bishop surface that is formally but not holomorphically equivalent to its quadric \cite{Gon04}. The construction reduces the obstruction to the linearization of a one-dimensional germ of holomorphic maps. The coefficients of the formal linearizer involve the denominators $\mu^n-\mu$, which can be very small. Gong chose a multiplier with extremely small divisors and picked the coefficients of an auxiliary holomorphic function $u$ one by one, forcing the formal linearizer coefficients to grow at least factorially along a subsequence. Gong's example has a further geometric property: it is not contained in any smooth Levi-flat real-analytic hypersurface \cite[Theorem 1.1]{Gon04}. His proof establishes existence but does not give an explicit example of such a surface. Further, the construction requires a specially chosen multiplier. Thus, the result shows that formal equivalence need not imply holomorphic equivalence. It leaves open which eigenvalues can occur in such counterexamples.

On the other hand, Gong had proved a convergence result in 1994 \cite{Gon94}. Specifically, suppose a real-analytic surface is formally equivalent to its hyperbolic Bishop quadric. If the associated multiplier satisfies a Diophantine condition, then this equivalence can be realized by a biholomorphic change of coordinates. Thus, formal and holomorphic equivalence agree under this arithmetic assumption.

A natural question is whether this Diophantine assumption can be weakened. The Brjuno condition \cite{Brj71} provides a weaker arithmetic hypothesis. Write $\mu=e^{2\pi i\alpha}$, where $\alpha\in\mathbb{R}\setminus\mathbb{Q}$. Let $q_n$ be the denominators of the continued-fraction convergents of $\alpha$. The Brjuno condition is $\sum_{n=0}^{\infty}\frac{\log q_{n+1}}{q_n}<\infty$. Gong and Stolovitch proved that a real-analytic surface formally equivalent to its hyperbolic Bishop quadric is holomorphically equivalent to it if the associated multiplier satisfies the Brjuno condition \cite[Theorem 6.7]{GS16}. This gives a convergence result for every Brjuno multiplier. Gong and Stolovitch conjectured that for every non-Brjuno unit multiplier $\mu$ that is not a root of unity, there exists a real-analytic hyperbolic Bishop surface with Moser--Webster invariant $\mu$ that is formally but not holomorphically equivalent to its quadric \cite[Conjecture 11.3]{GS26}.

Our main result proves this conjecture.

\begin{theorem}\label{thm:main}
Let $\mu\in\C$ satisfy $|\mu|=1$. Assume $\mu$ is neither a root of unity nor Brjuno. Let $\gamma=1/\sqrt{2+\mu+\mu^{-1}}>1/2$. Then there exists a germ of a real-analytic surface $M\subset\C^2$ with Moser--Webster invariant $\mu$ such that $M$ is formally equivalent, but not holomorphically equivalent, to its hyperbolic Bishop quadric $Q_\gamma$.
\end{theorem}

Together with the convergence theorem of Gong and Stolovitch, Theorem~\ref{thm:main} shows that the Brjuno condition is sharp.

We also show in Proposition~\ref{prop:levi-flat} that the surface $M$ constructed in Theorem~\ref{thm:main} cannot be contained in any smooth Levi-flat real-analytic hypersurface. A smooth Levi-flat real-analytic hypersurface is locally foliated by complex curves. Gong proved that, under the meromorphic eigenfunction condition satisfied by our construction, a surface contained in such a hypersurface must be holomorphically equivalent to its Bishop quadric \cite[Theorem 1.2]{Gon04}.

The main difference is that the multiplier can now be prescribed. Gong's earlier construction required a special choice of multiplier with extremely small divisors. Our result shows that this restriction is not necessary. For every non-Brjuno unit multiplier that is not a root of unity, we construct a real-analytic surface that is formally equivalent to its quadratic model but not holomorphically equivalent to it. Moreover, our counterexamples come from the special family introduced in \cite{Gon04} with an auxiliary function $u$ of degree one (see Section~\ref{sec:gong-family}). We do not need to choose the coefficients of $u$ one by one.

In fact, the small-divisor problem had already appeared in one-dimensional holomorphic dynamics. A germ can be formally linearizable without being holomorphically linearizable. In 1942, Siegel proved that holomorphic linearization holds when the multiplier satisfies a suitable Diophantine condition \cite{Sie42}. His result showed that arithmetic properties of the multiplier can determine whether a formal linearization converges. One-dimensional dynamics suggests that the Brjuno condition is the natural arithmetic threshold. Brjuno proved its sufficiency for holomorphic linearization, while Yoccoz proved its necessity for the quadratic polynomial $\mu z+z^2$ \cite{Brj71,Yoc95}.

We now outline the proof. The proof uses three main ingredients. First, Yoccoz's theorem shows that $P_\mu(z)=\mu z+z^2$ is not holomorphically linearizable for every non-Brjuno multiplier $\mu$. Second, we use a real one-parameter subfamily of Gong's construction. A non-linearizable germ in this family obstructs holomorphic linearization of the associated two-dimensional Moser--Webster map. Third, following P\'erez-Marco \cite{PM01}, we use the Bernstein--Walsh inequality to pass from bounds on a nonpolar compact set to bounds at each fixed complex parameter. We give the argument for our family in full.

To connect $P_\mu$ with Gong's construction, we define $a_t=\frac{1+it}{2}, b_t=\frac{1-it}{2}$, $A_t(z)=z+a_tz^2, B_t(z)=z+b_tz^2$, and
\begin{equation}\label{eq:Ft-intro}
F_t=A_t^{-1}\circ R_\mu\circ B_t,
\qquad R_\mu(z)=\mu z.
\end{equation}
The parameter has two roles: for real $t$, the germ $F_t$ is exactly Gong's one-dimensional germ, while at $t=i$ one has $F_i(z)=\mu(z+z^2)$, which is linearly conjugate to $P_\mu$. For nonreal $t$, we use $F_t$ only as a holomorphic complexification; we do not claim that it arises from a Moser--Webster surface.

The key fact is that the coefficients $h_n(t)$ of the normalized formal linearizer $H_t(z)=z+\sum_{n\ge2}h_n(t)z^n$ are polynomials in $t$ with $\deg h_n\le n-1$. The Bernstein--Walsh inequality shows that the set $\mathcal L_\mu$ of linearizable complex parameters is either all of $\C$ or a polar set. Yoccoz's theorem gives $i\notin\mathcal L_\mu$, so $\mathcal L_\mu$ is polar. A polar set cannot contain a nondegenerate real interval. Hence there is a real parameter for which Gong's germ is not holomorphically linearizable. This gives the one-dimensional obstruction needed for the Moser--Webster construction.

The paper is organized as follows. Section \ref{sec:inputs} collects the definitions and results used in the proof. Section \ref{sec:one-dimensional} proves the parameter dichotomy and applies it to the family \eqref{eq:Ft-intro} to establish Theorem \ref{thm:one-dimensional}. Section \ref{sec:surface} proves Theorem \ref{thm:main} and gives a consequence for four-dimensional real submanifolds in $\C^4$.

\section{Preliminaries}\label{sec:inputs}

\subsection{The Brjuno condition and Yoccoz's theorem}

Let $\alpha\in\R\setminus\mathbb Q$ and let $p_n/q_n$ be the convergents of its continued fraction. The Brjuno sum is $B(\alpha)=\sum_{n\ge0}\frac{\log q_{n+1}}{q_n}$.
We say that $\alpha$ is Brjuno if $B(\alpha)<\infty$ \cite{Brj71}. If $\mu=e^{2\pi i\alpha}$, we also say that $\mu$ is Brjuno.

Gong and Stolovitch use the following small-divisor condition \cite[Eq.~(10.14)]{GS26}:
$$
\sum_{k=1}^{\infty}2^{-k}
\log\frac{1}{\min_{1\le j\le 2^k}|\mu^j-1|}
<\infty.
$$
For irrational $\alpha$, this condition is equivalent to $B(\alpha)<\infty$; see Brjuno \cite[Chapter~II, \S5, Lemmas~12--13 and Example~7]{Brj71} and Ponce \cite[Section~2.1]{Pon09}.

Consider the quadratic polynomial $P_\mu(z)=\mu z+z^2$. We use the following consequence of Yoccoz's theorem \cite{Yoc95}.

\begin{theorem}\label{thm:yoccoz}
If $\mu=e^{2\pi i\alpha}$ with $\alpha\in\R\setminus\mathbb Q$ and $\alpha$ is not Brjuno, then the germ $P_\mu(z)=\mu z+z^2$ at the origin is not holomorphically linearizable.
\end{theorem}

Thus, for every non-Brjuno multiplier, Yoccoz supplies a concrete non-linearizable germ. At $t=i$, the germ $F_i$ is linearly conjugate to this quadratic. We use no quantitative estimate from Yoccoz's proof.

\subsection{Polar sets, logarithmic capacity, and real intervals}

We recall the potential-theoretic notions used in the parameter argument.

\begin{definition}\label{def:polar}
A set $E\subset\C$ is called polar if there exists a subharmonic function $U$ on $\C$, not identically $-\infty$, such that $U(z)=-\infty$ for every $z\in E$. Here a subharmonic function is an upper semicontinuous function $U:\C\to[-\infty,\infty)$ satisfying the submean inequality
$$
U(z_0)\le \frac{1}{2\pi}\int_0^{2\pi}U(z_0+re^{i\theta})\,d\theta
$$
for every $z_0\in\C$ and $r>0$.
\end{definition}

Let $K\subset\C$ be a nonempty compact set, and let $\mathcal P(K)$ be the set of Borel probability measures supported on $K$. Thus, $\nu\in\mathcal P(K)$ means that $\nu$ is a nonnegative Borel measure with $\nu(K)=1$ and $\nu(\C\setminus K)=0$. The logarithmic energy of $\nu\in\mathcal P(K)$ is
$$
I(\nu)=\iint_{K\times K}\log\frac1{|z-w|}\,d\nu(z)\,d\nu(w).
$$
The integrand is taken to be $+\infty$ when $z=w$. Since $K$ is compact, the integrand is bounded below on $K\times K$, so the integral is well defined, possibly with value $+\infty$. Define
$$
W(K)=\inf_{\nu\in\mathcal P(K)}I(\nu),
\qquad
\caplog(K)=e^{-W(K)},
$$
with $e^{-\infty}=0$. The number $\caplog(K)$ is the logarithmic capacity of $K$.

We use the following standard connection between the two definitions: a compact set $K$ is polar if and only if $\caplog(K)=0$. See \cite[Sections~3.2 and~5.1, pp.~55--57 and~127--128]{Ran95}; this equivalence is also stated explicitly in \cite[Section~4, p.~119]{BR18}. We also use that a countable union of polar sets is polar; see \cite[p.~459]{PM01}.

The following lemma will be used to find a real parameter for which linearization fails.

\begin{lemma}\label{lem:interval-capacity}
Every closed real interval of positive length has positive logarithmic capacity. More precisely, if $a<b$, then
$$
\caplog([a,b])\ge (b-a)e^{-3/2}>0.
$$
In particular, no polar subset of $\C$ can contain a nondegenerate real interval, and no polar subset can contain the whole real axis.
\end{lemma}

\begin{proof}
Fix $a<b$, put $L=b-a>0$ and $K=[a,b]$, and let $\nu$ be the uniform probability measure on $K$. Thus $d\nu(x)=dx/L$ for $a\le x\le b$, where $dx$ denotes one-dimensional Lebesgue measure, and $\nu(K)=L^{-1}\int_a^b dx=1$.

By definition, $I(\nu)=L^{-2}\int_a^b\int_a^b\log(1/|x-y|)\,dy\,dx$. The singularity at $x=y$ is integrable, since $\int_0^L|\log t|\,dt<\infty$. Put $x=a+s$ and $y=a+r$. Symmetry about the diagonal $s=r$ gives
$$
I(\nu)=\frac2{L^2}\int_0^L\int_0^s\log\frac1{s-r}\,dr\,ds.
$$
Now set $t=s-r$. The region of integration becomes $0\le t\le s\le L$. Changing the order of integration, we obtain
$$
I(\nu)=\frac2{L^2}\int_0^L(L-t)\log\frac1t\,dt.
$$
Integration by parts gives
$$
\int_0^L\log\frac1t\,dt=L-L\log L,\qquad
\int_0^L t\log\frac1t\,dt=\frac{L^2}{4}-\frac{L^2}{2}\log L.
$$
The antiderivatives are $-t\log t+t$ and $-t^2\log t/2+t^2/4$, respectively. The limits at $t=0$ exist because $t\log t\to0$ and $t^2\log t\to0$ as $t\downarrow0$. Substitution gives
$$
I(\nu)=\frac2{L^2}\left[L(L-L\log L)-\left(\frac{L^2}{4}-\frac{L^2}{2}\log L\right)\right]
=\frac32-\log L<\infty.
$$
Since $\nu\in\mathcal P(K)$, we have $W(K)\le I(\nu)$. Therefore
$$
\caplog(K)=e^{-W(K)}\ge e^{-I(\nu)}=Le^{-3/2}>0.
$$
By the equivalence between polarity and zero logarithmic capacity, $K$ is not polar. Every subset of a polar set is polar, so no polar set can contain a nondegenerate real interval. In particular, no polar set can contain the whole real axis.
\end{proof}

The exact capacity of an interval is also known:
$$
\caplog([a,b])=\frac{b-a}{4},
\qquad
\caplog([-1,1])=\frac12;
$$
see \cite[Example~1.11, p.~174]{Saff10}. We do not need these exact values below.

\subsection{The Bernstein--Walsh inequality}

A nondegenerate real interval has positive logarithmic capacity, despite having zero planar Lebesgue measure. The Bernstein--Walsh inequality explains why such a set is large enough to control a polynomial away from it.

Let $K\subset\C$ be a nonpolar compact set. For any polynomial $P$, write $\|P\|_K=\sup_{s\in K}|P(s)|$. We write $V_K$ for the one-variable Siciak extremal function, defined by
$$
V_K(t)=\sup\left(\{0\}\cup
\left\{
\frac{\log|Q(t)|}{\deg Q}:
Q\in\C[t],\ \deg Q\ge1,\ \|Q\|_K\le1
\right\}\right).
$$
The function $V_K$ is finite at every $t\in\C$ when $K$ is nonpolar. For any polynomial $P$ of degree at most $d$, the Bernstein--Walsh inequality gives
\begin{equation}\label{eq:BW}
|P(t)|\le \|P\|_K\exp\bigl(dV_K(t)\bigr),
\qquad t\in\C.
\end{equation}
This is the classical Bernstein--Walsh inequality; see \cite[p.~156]{Ran95} and \cite[Section~1.1, p.~459]{PM01}. It shows that an exponential bound on $\|P\|_K$ gives an exponential bound on $|P(t)|$ at any fixed complex parameter $t$, provided that the degree of $P$ is controlled.

As P\'erez-Marco notes, the use of polynomial parameter dependence in small-divisor problems goes back to Ilyashenko \cite{Ily79}; see \cite[Introduction, pp.~451--452]{PM01}. Our argument follows P\'erez-Marco \cite[Proof of Theorem~1, pp.~460--461]{PM01}. He combines this dependence with the Bernstein--Walsh inequality to show that either the formal linearizers converge for every parameter, or the set of convergent parameters is polar. We give the proof for our family for the convenience of the reader.

We apply \eqref{eq:BW} to the polynomial coefficients $h_n(t)$ of the formal linearizer, whose degrees are at most $n-1$. If the convergence set is nonpolar, we will find a nonpolar compact subset $K$ and an integer $N\ge1$ such that $\|h_n\|_K\le N^n$ for every $n\ge2$. Then \eqref{eq:BW} gives
$$
|h_n(t)|\le N^n\exp\bigl((n-1)V_K(t)\bigr)
\le\bigl(Ne^{V_K(t)}\bigr)^n
$$
for each fixed $t\in\C$. Since $V_K(t)<\infty$, this is an exponential bound in $n$. Thus the formal linearizer converges at $t$. The proof of the parameter dichotomy is given in Section \ref{sec:one-dimensional}.

\subsection{Gong's family of involutions}\label{sec:gong-family}

We recall Gong's construction \cite{Gon04}.

For a holomorphic germ $u(z)=\sum_{n\ge0}u_nz^n$, write
$$
\overline u(z)=\overline{u(\overline z)}
=\sum_{n\ge0}\overline{u_n}z^n.
$$
If $u(0)=1$, let $v$ be the unique holomorphic germ satisfying
\begin{equation}\label{eq:uv}
u(z)v(zu(z))=1.
\end{equation}

To see that $v$ exists and is unique, note that the map $z\mapsto zu(z)$ has derivative $1$ at $0$ and hence has a holomorphic inverse near $0$. This inverse has the form $z\mapsto zv(z)$, where $v$ is holomorphic and $v(0)=1$.
Equation \eqref{eq:uv} follows by composing the two maps.

Let $|\lambda|=1$ and put $\mu=\lambda^2$. Define $\rho(\xi,\eta)=(\overline\xi,\overline\eta)$ and $\tau_1^*(\xi,\eta)=(\lambda\eta,\overline\lambda\,\xi)$. With $u$ as above, set
\begin{equation}\label{eq:psi}
\psi(\xi,\eta)=(\xi u(\xi),\eta u(\xi)).
\end{equation}
Then $\psi^{-1}(\xi,\eta)=(\xi v(\xi),\eta v(\xi))$. Indeed, the identity $u(z)v(zu(z))=1$ gives $v(zu(z))=1/u(z)$.

Following Gong, define two holomorphic involutions and their composition by
\begin{equation}\label{eq:tau-pair}
\tau_1=\psi^{-1}\tau_1^*\psi,
\qquad
\tau_2=\rho\tau_1\rho,
\qquad
\sigma=\tau_1\tau_2.
\end{equation}
A direct calculation gives
\begin{equation}\label{eq:tau1}
\tau_1(\xi,\eta)
=\bigl(\lambda\eta\,p(\xi,\eta),
\overline\lambda\,\xi\,p(\xi,\eta)\bigr),
\end{equation}
where
\begin{equation}\label{eq:p}
p(\xi,\eta)=u(\xi)v(\lambda\eta u(\xi)).
\end{equation}
In particular, $p(\xi,0)=u(\xi)$ and $p(0,\eta)=v(\lambda\eta)$.

We use two properties of this construction.

\begin{enumerate}[label=(\roman*)]
\item Let $m(\xi,\eta)=\xi/\eta$. Since $\psi$ preserves this quotient, we have $m\circ\tau_1=\mu/m$ and $m\circ\tau_2=\mu^{-1}/m$. Consequently, $m\circ\sigma=\mu^2m$, so $m$ is a meromorphic eigenfunction of $\sigma$. If $\mu$ is not a root of unity, Proposition~3.1 in \cite{Gon04} implies that $\sigma$ is formally linearizable.

\item The coordinate axis $\{\eta=0\}$ is invariant under $\sigma$. On this axis,
\begin{equation}\label{eq:gong-axis}
\sigma(\xi,0)=(f(\xi),0),
\qquad
f(\xi)=\mu\xi\,\overline u(\xi)
 v\!\left(\mu\xi\,\overline u(\xi)\right).
\end{equation}
If $\sigma$ is holomorphically linearizable, then $f$ is also holomorphically linearizable. We will give a direct calculation of $\sigma(\xi,0)$ and prove the implication about linearizability in Section \ref{sec:surface}.
\end{enumerate}

Thus, the remaining task is to choose a holomorphic germ $u$ for which the one-dimensional germ $f$ in \eqref{eq:gong-axis} is not holomorphically linearizable.

\subsection{The Moser--Webster realization}

We use the Moser--Webster correspondence \cite{MW83} in the form recalled by Gong \cite[Section~2]{Gon04} and Gong--Stolovitch \cite[Theorem~2.2]{GS26}. Choose $|\lambda|=1$ with $\Re\lambda>0$, and assume that $\mu=\lambda^2$ is not a root of unity. Then $\gamma=(\lambda+\lambda^{-1})^{-1}>1/2$. The linear parts of $\tau_1$ and $\tau_2$ are $\tau_1^*(\xi,\eta)=(\lambda\eta,\overline\lambda\xi)$ and $\tau_2^*(\xi,\eta)=(\overline\lambda\eta,\lambda\xi)$, respectively. These involutions satisfy $\tau_2=\rho\tau_1\rho$. The Moser--Webster realization theorem therefore gives a germ of a real-analytic Bishop surface with Bishop invariant $\gamma$.

Under this correspondence, two surfaces are formally (respectively, holomorphically) equivalent if and only if their associated involution triples are simultaneously conjugate by a formal (respectively, holomorphic) coordinate change commuting with $\rho$. Since $\mu$ is not a root of unity, formal linearization of $\sigma=\tau_1\tau_2$ also implies simultaneous formal linearization of the involution triple, by the Moser--Webster formal normal form \cite[Theorem~2.3]{GS26}. We apply this correspondence in Section \ref{sec:surface} to prove Theorem \ref{thm:main}.

\section{The one-dimensional obstruction}\label{sec:one-dimensional}

\subsection{A parameter dichotomy}\label{sec:parameter}

The argument in this subsection is essentially the argument of P\'erez-Marco \cite[Proof of Theorem~1, pp.~460--461]{PM01}. We give the proof for our setting for the convenience of the reader.

\begin{proposition}\label{prop:parameter}
Let $\mu\in\C\setminus\{0\}$, and assume that $\mu$ is not a root of unity. For each $t\in\C$, let $F_t(z)=\mu z+\sum_{n\ge2}c_n(t)z^n$ be a convergent holomorphic germ at the origin. Assume that every $c_n(t)$ is a polynomial in $t$ with $\deg c_n\le n-1$.

Let $H_t(z)=z+\sum_{n\ge2}h_n(t)z^n$ be the unique normalized formal solution of $H_t(F_t(z))=\mu H_t(z)$. Then:
\begin{enumerate}[label=(\alph*)]
\item Every $h_n(t)$ is a polynomial satisfying $\deg h_n\le n-1$.
\item The set $\mathcal L$ of parameters $t\in\C$ for which $F_t$ is holomorphically linearizable is either $\C$ or a polar set.
\end{enumerate}
\end{proposition}

\begin{proof}
We first prove part (a). Put $c_1(t)=\mu$. For $k\ge2$, write
$$
F_t(z)^k=\sum_{n\ge k}C_{n,k}(t)z^n.
$$
Each $C_{n,k}(t)$ is a sum of products $c_{j_1}(t)\cdots c_{j_k}(t)$, where $j_1,\ldots,j_k\ge1$ and $j_1+\cdots+j_k=n$. Since $\deg c_j\le j-1$, each product has degree at most $(j_1-1)+\cdots+(j_k-1)=n-k$. Hence $\deg C_{n,k}\le n-k$.

Comparing the coefficients of $z^n$ in $H_t(F_t(z))=\mu H_t(z)$ gives
\begin{equation}\label{eq:hn-recursion}
(\mu-\mu^n)h_n(t)=c_n(t)+\sum_{k=2}^{n-1}h_k(t)C_{n,k}(t),\qquad n\ge2.
\end{equation}
Since $\mu$ is not a root of unity, $\mu-\mu^n\ne0$ for every $n\ge2$. Thus the equation determines each $h_n$ uniquely.

For $n=2$, the sum is empty, and $\deg h_2\le1$. Assume that $\deg h_k\le k-1$ for $2\le k<n$. Then every product $h_kC_{n,k}$ in \eqref{eq:hn-recursion} has degree at most $(k-1)+(n-k)=n-1$. The same bound holds for $c_n$. Hence $h_n$ is a polynomial with $\deg h_n\le n-1$. This proves part (a).

For part (b), $F_t$ is holomorphically linearizable if and only if $H_t$ converges. Indeed, a convergent $H_t$ has $H_t'(0)=1$ and gives a local holomorphic conjugacy. Conversely, any holomorphic linearizer can be normalized to have derivative $1$ at the origin. By uniqueness, its Taylor series is $H_t$. Thus
\begin{equation}\label{eq:lin-iff-conv}
F_t\text{ is holomorphically linearizable}
\quad\Longleftrightarrow\quad
H_t\text{ converges near }0.
\end{equation}

Let $\mathcal L=\{t\in\C:H_t\text{ converges near }0\}$. For each integer $N\ge1$, define
$$
K_N=\{t\in\C:|t|\le N,\ |h_n(t)|\le N^n\text{ for all }n\ge2\}.
$$
Each $K_N$ is compact, since it is bounded and every $h_n$ is continuous. Moreover, $\mathcal L=\bigcup_{N\ge1}K_N$. Indeed, if $t\in K_N$, the series $H_t$ converges for $|z|<1/N$. Conversely, if $H_t$ converges, then $|h_n(t)|\le M^n$ for some $M>0$ and all $n\ge2$. Taking an integer $N\ge\max\{1,|t|,M\}$ gives $t\in K_N$.

Suppose that $\mathcal L$ is not polar. Since a countable union of polar sets is polar, at least one $K_N$ is nonpolar. Fix such an $N$ and put $K=K_N$. By definition, $\|h_n\|_K=\sup_{s\in K}|h_n(s)|\le N^n$ for every $n\ge2$. Part (a) gives $\deg h_n\le n-1$. The Bernstein--Walsh inequality \eqref{eq:BW} therefore gives
$$
|h_n(t)|\le N^n\exp\bigl((n-1)V_K(t)\bigr)
\le\bigl(Ne^{V_K(t)}\bigr)^n,
\qquad t\in\C.
$$
Since $V_K(t)<\infty$ for each $t\in\C$, the series $H_t$ converges for every parameter. Thus $\mathcal L=\C$, and the dichotomy follows.
\end{proof}

\subsection{The special parameter family}\label{sec:family}

We now construct a one-parameter family that connects $P_\mu$ with Gong's one-dimensional germ. The complex parameter allows us to apply Proposition \ref{prop:parameter}, while real parameters will give the germs needed for the geometric construction.

Fix $\mu\in\C$ with $|\mu|=1$, and assume that $\mu$ is not a root of unity. For $t\in\C$, set $a_t=(1+it)/2$ and $b_t=(1-it)/2$. Define $A_t(z)=z+a_tz^2$ and $B_t(z)=z+b_tz^2$. Since $A_t'(0)=1$, the inverse $A_t^{-1}$ is a convergent holomorphic germ for every fixed $t$. We define
\begin{equation}\label{eq:F-special}
F_t(z)=A_t^{-1}(\mu B_t(z)).
\end{equation}
Thus, for every $t\in\C$, the map $F_t$ is a convergent holomorphic germ at the origin, with $F_t(0)=0$ and $F_t'(0)=\mu$.

To apply the parameter dichotomy, we first check how the coefficients of $F_t$ depend on $t$. The defining equation gives
\begin{equation}\label{eq:family-identity}
F_t(z)+a_tF_t(z)^2=\mu z+\mu b_tz^2.
\end{equation}
Write $F_t(z)=\sum_{n\ge1}c_n(t)z^n$, where $c_1(t)=\mu$. Comparing coefficients in \eqref{eq:family-identity}, we obtain $c_2(t)=\mu b_t-\mu^2a_t$ and, for every $n\ge3$,
\begin{equation}\label{eq:cn-recursion}
c_n(t)=-a_t\sum_{j=1}^{n-1}c_j(t)c_{n-j}(t).
\end{equation}

The next lemma verifies the degree bound required by Proposition \ref{prop:parameter}. It follows directly from this recursion.

\begin{lemma}\label{lem:cn-degree}
For every $n\ge2$, the coefficient $c_n(t)$ is a polynomial in $t$ satisfying $\deg c_n\le n-1$.
\end{lemma}

\begin{proof}
Since $a_t$ and $b_t$ are affine in $t$, the assertion holds for $n=2$. Suppose that $\deg c_j\le j-1$ for every $2\le j<n$. Since $c_1(t)=\mu$ is constant, each product $c_j(t)c_{n-j}(t)$ in \eqref{eq:cn-recursion} has degree at most $(j-1)+(n-j-1)=n-2$. Multiplication by $a_t$ increases the degree by at most one. Hence $\deg c_n\le n-1$, completing the induction.
\end{proof}

We next identify two useful parameters. At $t=0$, we have $a_0=b_0=1/2$ and hence $A_0=B_0$. It follows that $F_0=A_0^{-1}\circ R_\mu\circ A_0$, where $R_\mu(z)=\mu z$. Thus $F_0$ is holomorphically linearizable.

At $t=i$, we have $a_i=0$ and $b_i=1$. Therefore
\begin{equation}\label{eq:ti}
F_i(z)=\mu(z+z^2).
\end{equation}
Let $P_\mu(z)=\mu z+z^2$ and $L(z)=\mu z$. A direct calculation gives $F_i=L^{-1}\circ P_\mu\circ L$. Thus $F_i$ is linearly conjugate to $P_\mu$.

We now relate \eqref{eq:F-special} to Gong's construction. For real $t$, set $u_t(z)=1+a_tz$, and let $v_t$ be determined by \eqref{eq:uv}. Since $b_t=\overline{a_t}$, we have $A_t(z)=zu_t(z)$ and $B_t(z)=z\overline{u_t}(z)$. The defining relation between $u_t$ and $v_t$ also gives $A_t^{-1}(w)=wv_t(w)$. Consequently,
\begin{equation}\label{eq:real-slice}
F_t(z)=\mu z\,\overline{u_t}(z)
v_t\!\left(\mu z\,\overline{u_t}(z)\right)
=f_t(z),\qquad t\in\R.
\end{equation}
This is exactly the one-dimensional germ in Gong's construction \cite{Gon04}, with the auxiliary function $u=u_t$.

\noindent\emph{Remark.} For nonreal $t$, the relation $b_t=\overline{a_t}$ no longer holds. We use these parameters only as a holomorphic complexification of the real family, and do not claim that the corresponding germs arise from Moser--Webster surfaces.

We can now prove the one-dimensional result needed for the main theorem. It shows that the real parameters for which Gong's germ is holomorphically linearizable form a polar set. In particular, every real interval of positive length contains a parameter for which holomorphic linearization fails.

\begin{theorem}\label{thm:one-dimensional}
Let $\mu\in\C$ satisfy $|\mu|=1$, and assume that $\mu$ is neither a root of unity nor Brjuno. Let $E_\mu$ be the set of real parameters $t$ for which $f_t$ is holomorphically linearizable. Then $E_\mu$, viewed as a subset of $\C$, is polar. Consequently, every real interval of positive length contains a parameter $t$ for which $f_t$ is formally but not holomorphically linearizable.
\end{theorem}

\begin{proof}
By Lemma \ref{lem:cn-degree}, the coefficients of $F_t$ satisfy the degree assumptions of Proposition \ref{prop:parameter}. Moreover, $F_t$ is a convergent holomorphic germ for every $t\in\C$, with the same nonresonant multiplier $\mu$. Thus, the proposition applies. Let $\mathcal L_\mu=\{t\in\C:F_t\text{ is holomorphically linearizable}\}$. Proposition \ref{prop:parameter} shows that $\mathcal L_\mu=\C$ or $\mathcal L_\mu$ is polar.

Since $\mu$ is not Brjuno, Theorem \ref{thm:yoccoz} shows that $P_\mu$ is not holomorphically linearizable. By \eqref{eq:ti}, the germ $F_i$ is linearly conjugate to $P_\mu$. Hence $i\notin\mathcal L_\mu$. The first alternative is therefore impossible, and $\mathcal L_\mu$ must be polar. For real $t$, equation \eqref{eq:real-slice} gives $F_t=f_t$. Thus $E_\mu=\mathcal L_\mu\cap\R$, so $E_\mu$ is also polar as a subset of $\C$. Let $I$ be any real interval of positive length. It contains a closed interval $[a,b]$ with $a<b$, which is nonpolar by Lemma \ref{lem:interval-capacity}. Thus $I$ cannot be contained in $E_\mu$. Hence there exists $t\in I\setminus E_\mu$, for which $f_t$ is not holomorphically linearizable.

Finally, $f_t'(0)=\mu$, and $\mu$ is not a root of unity. The formal linearization equation therefore has a unique normalized formal solution, as follows from the coefficient recursion \eqref{eq:hn-recursion}. Thus $f_t$ is formally linearizable but not holomorphically linearizable. This proves the theorem.
\end{proof}

\section{Proof of the main theorem}\label{sec:surface}

We now use the one-dimensional result to construct the required Bishop surface. The main point is to show that failure of linearization on an invariant axis prevents holomorphic linearization of the full two-dimensional map. We first record this fact.

\begin{lemma}\label{lem:axis-obstruction}
Let $\mu\in\C\setminus\{0\}$, and assume that $\mu$ is not a root of unity. Let $S$ be a holomorphic germ at the origin of $\C^2$ with linear part $S^*(\xi,\eta)=(\mu\xi,\mu^{-1}\eta)$. Assume that both coordinate axes are invariant under $S$. If $S$ is holomorphically linearizable, then its restriction to each coordinate axis is holomorphically linearizable.
\end{lemma}

\begin{proof}
We first show that the coordinate axes are the only smooth invariant holomorphic curves through the origin for $S^*$. The tangent line of such a curve must be invariant under $DS^*(0)$. Since $\mu\ne\mu^{-1}$, it must be one of the two coordinate axes.

A curve tangent to the $\xi$-axis can be written as $\eta=q(\xi)$, where $q(0)=q'(0)=0$. Invariance under $S^*$ gives $q(\mu\xi)=\mu^{-1}q(\xi)$. If $q$ is not identically zero, write $q(\xi)=c\xi^r+O(\xi^{r+1})$, where $r\ge2$ and $c\ne0$. Comparing the first nonzero terms gives $\mu^{r+1}=1$, a contradiction. Hence $q=0$. The same argument applies to a curve tangent to the $\eta$-axis. Thus the coordinate axes are the only such invariant curves.

Suppose now that $S$ is holomorphically linearizable. Choose a local biholomorphism $\Phi$ satisfying $\Phi^{-1}S\Phi=S^*$. Since $DS(0)=DS^*(0)$, differentiation at the origin shows that $D\Phi(0)$ commutes with $DS^*(0)$. Since $DS^*(0)$ has distinct eigenvalues, $D\Phi(0)$ is diagonal. The inverse image under $\Phi$ of either invariant axis of $S$ is therefore a smooth invariant curve for $S^*$ tangent to the same coordinate axis. By the uniqueness proved above, $\Phi$ preserves both axes.

In particular, $\Phi(\xi,0)=(g(\xi),0)$ for some one-variable biholomorphic germ $g$. Restricting the conjugacy identity to this axis gives $g^{-1}\circ(S|_{\eta=0})\circ g(\xi)=\mu\xi$. Thus the restriction of $S$ to $\eta=0$ is holomorphically linearizable. The other axis is treated in the same way.
\end{proof}

We now choose the parameter and apply Gong's construction. The resulting involution triple will be formally equivalent to the linear triple, while the lemma will rule out holomorphic equivalence.

\begin{proof}[Proof of Theorem \ref{thm:main}]
Fix $\mu$ as in the theorem. Choose a square root $\lambda$ of $\mu$ with $\Ree\lambda>0$. This is possible because $\mu\ne-1$. Since $\mu$ is not a root of unity, we have $0<\Ree\lambda<1$. Set $\gamma=1/(2\Ree\lambda)$, so $\gamma>1/2$. As $|\lambda|=1$, we have $\lambda+\lambda^{-1}=1/\gamma$. Hence $\gamma\lambda^2-\lambda+\gamma=0$ and $\mu+\mu^{-1}=\gamma^{-2}-2$. Thus $\gamma$ is the hyperbolic Bishop invariant corresponding to the prescribed multiplier $\mu$.

By Theorem \ref{thm:one-dimensional}, we may choose a real parameter $t\notin E_\mu$. Set $u(z)=1+(1+it)z/2$, and let $v$ be determined by \eqref{eq:uv}. Define $\psi,\tau_1,\tau_2$ and $\sigma=\tau_1\tau_2$ as in \eqref{eq:psi} and \eqref{eq:tau-pair}.

We first establish formal equivalence. Since $\psi$ is tangent to the identity, the linear parts of the involutions are
\begin{equation*}
\tau_1^*(\xi,\eta)
=(\lambda\eta,\overline\lambda\,\xi),
\qquad
\tau_2^*(\xi,\eta)
=(\overline\lambda\,\eta,\lambda\xi).
\end{equation*}
Both maps are holomorphic involutions, and $\tau_2=\rho\tau_1\rho$. Their composition has linear part $\sigma^*(\xi,\eta)=(\mu\xi,\mu^{-1}\eta)$.

Consider the meromorphic function $m(\xi,\eta)=\xi/\eta$. The maps $\psi$ and $\psi^{-1}$ preserve this quotient. As shown in Section \ref{sec:inputs}, we have $m\circ\tau_1=\mu/m$ and $m\circ\tau_2=\mu^{-1}/m$. Consequently, $m\circ\sigma=\mu^2m$.

The numerator and denominator of $m$ are $\xi$ and $\eta$, respectively. Since $\mu=\lambda^2$ is not a root of unity, neither is $\lambda$. Thus Gong's Proposition 3.1 \cite[pp.~87--88]{Gon04} applies, and $\sigma$ is formally conjugate to $\sigma^*$.

We also need formal equivalence of the full involution triples. By the Moser--Webster formal normal form \cite{MW83}, as stated in \cite[Theorem~2.3]{GS26}, there is a unique normalized formal change of coordinates $\Phi$, tangent to the identity and commuting with $\rho$, in which the two involutions have scalar factors $\Lambda(\xi\eta)$ and $\Lambda(\xi\eta)^{-1}$, with $\Lambda(0)=\lambda$. Their composition has the form
\begin{equation*}
\widehat\sigma(\xi,\eta)
=\bigl(\Lambda(\xi\eta)^2\xi,
\Lambda(\xi\eta)^{-2}\eta\bigr).
\end{equation*}
The normalization of $\Phi$ sets the coefficients of $\xi^{j+1}\eta^j$ in its first component and of $\xi^j\eta^{j+1}$ in its second component equal to zero for every $j\ge1$. These are the same conditions used in Gong's Proposition 3.1 \cite[pp.~87--88]{Gon04}. In its proof, Gong shows that the unique normalized form of $\sigma$ is linear when the meromorphic eigenfunction exists. Therefore $\widehat\sigma=\sigma^*$ and $\Lambda(\xi\eta)^2=\mu$. Since $\Lambda(0)=\lambda$, we have $\Lambda=\lambda$. Thus the normalized involutions are the linear pair $\tau_1^*,\tau_2^*$. As $\Phi$ commutes with $\rho$, the entire triple $(\tau_1,\tau_2,\rho)$ is formally equivalent to $(\tau_1^*,\tau_2^*,\rho)$.

We next show that $\sigma$ is not holomorphically linearizable. The map $\rho$ preserves both coordinate axes, while $\tau_1$ interchanges them by \eqref{eq:tau1}. Hence $\sigma$ preserves each axis.

To compute its restriction to $\eta=0$, recall that $p(\xi,0)=u(\xi)$ and $p(0,\eta)=v(\lambda\eta)$. The formulas for $\tau_1$ and $\tau_2=\rho\tau_1\rho$ give
$\tau_2(\xi,0)=(0,\lambda\xi\overline u(\xi))$.
Applying $\tau_1$ once more, we obtain
\begin{equation*}
\sigma(\xi,0)
=\bigl(\mu\xi\overline u(\xi)
v(\mu\xi\overline u(\xi)),0\bigr)
=(f_t(\xi),0).
\end{equation*}
This is the real-slice formula \eqref{eq:real-slice}. Since $t\notin E_\mu$, the germ $f_t$ is not holomorphically linearizable.

If $\sigma$ were holomorphically linearizable, Lemma \ref{lem:axis-obstruction} would imply that its restriction $f_t$ to the invariant axis $\eta=0$ is holomorphically linearizable. This is a contradiction. Therefore $\sigma$ is not holomorphically linearizable.

We now return to the surface. The involutions $\tau_1,\tau_2$ have the linear parts displayed above and satisfy $\tau_2=\rho\tau_1\rho$. Since $\gamma>1/2$ and $\mu$ is not a root of unity, the Moser--Webster realization theorem \cite[Theorem~2.2]{GS26} gives a germ of a real-analytic surface $M\subset\C^2$ with Bishop invariant $\gamma$ and associated involution triple $(\tau_1,\tau_2,\rho)$.

Since the linear part of $\sigma$ has eigenvalues $\mu$ and $\mu^{-1}$, the surface has the prescribed Moser--Webster invariant. The formal equivalence of the involution triples gives a formal equivalence between $M$ and $Q_\gamma$.

If $M$ were holomorphically equivalent to $Q_\gamma$, the Moser--Webster correspondence would give a holomorphic conjugacy between their involution triples. In particular, $\sigma$ would be holomorphically linearizable. We have already ruled this out. Thus $M$ is formally but not holomorphically equivalent to $Q_\gamma$.
\end{proof}

\begin{proposition}\label{prop:levi-flat}
The surface $M$ constructed in the proof of Theorem~\ref{thm:main} is not contained in any smooth Levi-flat real-analytic hypersurface.
\end{proposition}

\begin{proof}
We have $\gamma>1/2$, and $\lambda$ is not a root of unity. Also, $\sigma=\tau_1\tau_2$ has the nonconstant meromorphic eigenfunction $m(\xi,\eta)=\xi/\eta$, with $m\circ\sigma=\mu^2m$. Its numerator and denominator are the holomorphic functions $\xi$ and $\eta$. Thus the hypotheses of Gong's theorem \cite[Theorem~1.2, p.~84]{Gon04} are satisfied.

If $M$ were contained in a smooth Levi-flat real-analytic hypersurface, that theorem would imply that $M$ is holomorphically equivalent to $Q_\gamma$, contrary to what we have proved. Therefore $M$ is not contained in any smooth Levi-flat real-analytic hypersurface.
\end{proof}

\begin{corollary}\label{cor:product}
Let $\gamma_1$ be the Bishop invariant in Theorem~\ref{thm:main}, and let $\gamma_2>1/2$. There exists a germ of a real-analytic submanifold $\mathcal M\subset\C^4$ of real dimension $4$ that is formally, but not holomorphically, equivalent to $Q_{\gamma_1}\times Q_{\gamma_2}$.
\end{corollary}

\begin{proof}
Let $M\subset\C^2$ be the surface constructed in Theorem~\ref{thm:main}, and set $\mathcal M=M\times Q_{\gamma_2}$. The formal equivalence $M\sim Q_{\gamma_1}$, extended by the identity on the second factor, gives a formal equivalence $\mathcal M\sim Q_{\gamma_1}\times Q_{\gamma_2}$.

Suppose that a local biholomorphism $H=(F,G)$ maps $\mathcal M$ onto $Q_{\gamma_1}\times Q_{\gamma_2}$. Write the coordinates of $\C^4$ as $(z,w)\in\C^2\times\C^2$. The second-order truncation of a formal equivalence between $M$ and $Q_{\gamma_1}$ defines a polynomial map with invertible derivative. This map is a local biholomorphism and brings the quadratic part of $M$ to that of $Q_{\gamma_1}$. We use these coordinates in the first factor, so that
$$
M:\quad w_1=z_1\overline z_1+\gamma_1(z_1^2+\overline z_1^{\,2})+O(|z_1|^3).
$$
Since both tangent spaces at the origin are $\{w=0\}$, the derivative of $H$ has the form
$$
DH(0)=\begin{pmatrix}A&B\\0&D\end{pmatrix},
\qquad A,D\in\operatorname{GL}_2(\C).
$$
Comparing the terms of bidegree $(1,1)$ in the two defining equations gives, for $j=1,2$,
$$
\sum_{k=1}^2 D_{jk}z_k\overline z_k
=\left|\sum_{k=1}^2 A_{jk}z_k\right|^2.
$$
Thus $D_{jk}=|A_{jk}|^2$ and $A_{j1}\overline{A_{j2}}=0$. Since $A$ is invertible, each row and each column has exactly one nonzero entry. Choose $j_0\in\{1,2\}$ so that $a:=A_{j_0 1}\ne0$. Then $D_{j_0 1}=|a|^2$.

Comparing the coefficients of $\overline z_1^{\,2}$ in the $j_0$-th defining equation gives $D_{j_0 1}\gamma_1=\gamma_{j_0}\overline a^{\,2}$. Since $D_{j_0 1}=|a|^2$ and $\gamma_1,\gamma_{j_0}>0$, taking absolute values shows that $\gamma_{j_0}=\gamma_1$.

Let $\iota:\C^2\to\C^4$ insert a point into the first factor, and let $\pi_{j_0}:\C^4\to\C^2$ project onto the $j_0$-th factor. Set $h=\pi_{j_0}\circ H\circ\iota$. Then $\det Dh(0)=aD_{j_0 1}=a|a|^2\ne0$, so $h$ is a local biholomorphism. Since $H(\mathcal M)\subset Q_{\gamma_1}\times Q_{\gamma_2}$, we have $h(M)\subset Q_{\gamma_{j_0}}=Q_{\gamma_1}$. Both are germs of real surfaces of dimension $2$, so this inclusion is an equality of germs. Thus $h$ is a holomorphic equivalence between $M$ and $Q_{\gamma_1}$, contradicting Theorem~\ref{thm:main}.
\end{proof}

Our optimal result is sharp for real surfaces in $\C^2$. However, we don't know if the Brjuno condition is sharp for higher dimensions.

\section*{Acknowledgments}

I am deeply grateful to Professor Xianghong Gong for his generous and sustained guidance throughout this project. Our frequent discussions played an essential role in the development of this work, and many important ideas emerged from them. This paper would not have been possible without his guidance and support.

\section*{Statements and Declarations}

The author has no relevant financial or non-financial interests to disclose.

\noindent\textbf{Data availability.}
No datasets were generated or analyzed during the current study.

\end{document}